**Working Title:** Math for America Los Angeles: Lessons Learned from 10 Years of Teacher-Centered Leadership Development
**Authors:**
1) Darryl H. Yong, Professor of Mathematics, Harvey Mudd College, 301 Platt Blvd, Claremont, CA 91711, 909-607-2844, dyong@hmc.edu
2) Pam Mason, Executive Director, Math for America Los Angeles, 1150 S Olive Street, Suite 2100, Los Angeles, CA 90015, 213-821-4362, pammason@mfala.org
3) Rebecca M. Eddy, President, Cobblestone Applied Research and Evaluation, Inc., 2120 Foothill Blvd, Suite 202, La Verne, CA 91750, 800-971-3891, rebecca.eddy@cobblestoneeval.com
4) Ashley Alchehayed, Research Assistant, Cobblestone Applied Research and Evaluation, Inc., 2120 Foothill Blvd, Suite 202, La Verne, CA 91750, 800-971-3891, ashley.alchehayed@cobblestoneeval.com
5) David Mendelsohn, Statistics Consultant, Cobblestone Applied Research and Evaluation, Inc., 2120 Foothill Blvd, Suite 202, La Verne, CA 91750, 800-971-3891, david.mendelsohn@cobblestoneeval.com

**Abstract:** This article describes what can happen when sustained effort and resources are devoted to creating a teacher professional support and development organization that puts teachers' needs first. Over the last ten years, Math for America Los Angeles has supported 179 secondary school mathematics and computer science teachers in developing their leadership through multi-year teaching fellowships. We share some of its history, design principles, and effects on teachers and students.

# 1. Introduction

Here is a paradox: though the field of teacher leadership has amassed wisdom about what effective teacher professional learning should look like, robust evidence of that effectiveness remains elusive (Garet et al., 2011; Hill, 2015; Jacob & McGovern, 2015). Yet, "one constant finding in the research literature is that notable improvements in education almost never take place *in the absence* of [teacher] professional development" (Guskey, 2000, p. 4).

In this article, we describe what can happen when sustained effort and resources are devoted to creating a "gold-standard" teacher professional support and development organization that puts teachers' needs first. Taking teachers' needs seriously involves listening to what teachers say they need and want for their professional learning, respecting them as professionals, and caring for them as individuals. Creating a high-quality professional development organization requires attending to that amassed wisdom about what effective teacher professional development looks like, investing sustained effort and resources in teachers over time, while attending carefully to program evaluation and assessment.

Math for America Los Angeles (MfA LA) is a non-profit organization founded in 2007 by the University of Southern California, Claremont Graduate University, and Harvey Mudd College that aims to raise student achievement by improving secondary school mathematics and computer science (CS) teacher recruitment, retention, and quality in the greater Los Angeles area. We aspire to achieve a critical mass of highly-effective mathematics and CS teacher leaders in the greater Los Angeles area so as to improve student learning and bring about more equity and justice. We strive to nurture teacher leaders because of the robust evidence that teachers have a significant effect on student learning (Hattie, 2009; Haycock, 1998; McCaffrey, J.R. Lockwood, Koretz, & Hamilton, 2003; Nye, Konstantopoulos, & Hedges, 2004; OECD, 2005; Timperley & Alton-Lee, 2008).

Why does MfA LA focus on mathematics and CS teachers? Mathematics is crucial to the preparation of future scientists, engineers, and creative thinkers and is a critical gate-keeper for many students pursuing STEM fields. At the same time, shifts in technology and employment trends have established the importance of CS as a basic skill necessary for full participation in the modern workforce. Access to CS courses in schools correlates highly with the percentage of Caucasian families and socio-economic status of families in the neighborhood (Margolis, 2010). Therefore, unequal access to quality math and CS education is a major barrier to broader participation in STEM. Moreover, high poverty, high minority, and urban public schools have the highest rates of mathematics and science teacher turnover and shortages and schools have become more segregated by race and income over the last two decades (Ingersoll & May, 2010; Owens, Reardon, & Jencks, 2016).

In Section 2, we explain how MfA LA builds on other teacher professional development efforts and reflects long-standing understandings about effective teacher professional development. Next, we describe MfA LA's program design and evidence of impact in Sections 3 and 4. Finally, we close with some lessons about teacher professional development that we have learned.

## 2. Building on the Other Teacher Professional Development Efforts

Two of us (Pam Mason and Darryl Yong) have been the primary architects of MfA LA's professional development efforts. Before we describe specific details about those efforts in Section 3, we first reflect on prior experiences that have affected how we think about teacher professional development. We do this to acknowledge and honor the efforts of others who have indirectly contributed to this work and to show how MfA LA relates to and benefits from the accumulated wisdom of prior and contemporary efforts. Finally, we also connect these ideas to research on the characteristics of successful teacher professional development programs.

### Teacher Leadership Opportunities and National Board Certification

Throughout my career as a middle school teacher, I (Pam) had many opportunities to develop my leadership skills. After the Los Angeles Unified School District (LAUSD) teachers' strike in 1989, negotiations with the district included restructuring our schools and implementing school site decision-making. The new contract in 1991 empowered teachers to be integral partners in school reform and management. Our school (Patrick Henry Middle School) was selected as one of the first eight Los Angeles Education Alliance for Restructuring Now (LEARN) schools. LEARN was an organization started by a small group of business executives and local activists that attempted to restructure schools to have greater local control. At that time, I was math department chair, math coach, the teacher's union (United Teachers of Los Angeles) chapter chair, union representative for our local area, and the lead teacher with my principal co-leading our school under LEARN. I was then asked by the union to be its mathematics educator representative on all professional development and curriculum committees and at district mathematics meetings. That was the first time in my career that teachers at each school site and district-wide were informing the decisions that were being made regarding professional development programming.

Around the same time, my school was selected to be part of Los Angeles Systemic Initiative, a federally-funded program designed to improve mathematics and science education. My classroom became the meeting place for math teachers in our area to discuss, reflect, work together to implement a new integrated mathematics curriculum. It was helpful that I was teaching and coaching at that time because I could model what worked in my own teaching. This teacher-initiated style of professional learning that avoids dogmatic insistence on one pedagogical approach reflects strongly in MfA LA's design.

Another formative experience for me was my involvement in a pilot program for the National Boards for Early Adolescence and Adolescence and Young Adulthood Mathematics in 1997. The process of becoming a National Board Certified Teacher was one of the best professional development experiences I have had. That experience is the reason why National Board certification is an integral part of MfA LA's Master Teacher Fellowship Program.

## Mathematics Institutes and Teacher Circles

In the late 1980s and early 1990s, the National Science Foundation Divisions of Mathematical Sciences and Elementary, Secondary and Informal Education promoted the idea of integrated mathematics institutes. A number of Regional Geometry Institutes resulted from those efforts. Herb Clemens, Dan Freed, and Karen Uhlenbeck created one of these institutes that became the IAS/Park City Mathematics Institute (PCMI). This intensive summer program brought mathematics researchers and educators together to learn amongst each other. From its start, PCMI included a program for school teachers under the direction of Jim King. Currently, that program is called the Teacher Leadership Program. It has roots that can be traced back to the Ross Program at Ohio State University and PROMYS at Boston University.

I (Darryl) first participated in a program for undergraduate mathematics faculty at PCMI in 2003. There, I met Peg Cagle, a middle school math teacher, who had started a mathematics teacher circle in Los Angeles and was looking for an undergraduate mathematics faculty to co-organize it. Math teacher circles were much less common at that time than they are now; it was unusual for a group of teachers from different districts to gather together on a Saturday to do mathematics together. My involvement in PCMI's Teacher Leadership Program and that Los Angeles mathematics teacher circle gave me a lasting, profound respect for the work that teachers do and helped me to understand the value of helping mathematics teachers develop their own mathematical identities. It was also through this teacher circle that I got to know Pam Mason.

## Brief History of Math for America Los Angeles (MfA LA)

Founded in January 2004 by Jim Simons, Math for America is a nonprofit organization whose mission is to promote recruitment and retention of high quality mathematics teachers in New York City secondary schools. Math for America spawned autonomous but affiliated organizations in San Diego, Los Angeles, Berkeley, Utah, Boston, and the District of Columbia, each with its own approach to teacher recruitment, certification, and professional development. Math for America Los Angeles was established in 2007 by Maria Klawe (president of Harvey Mudd College), Karen Gallagher (dean of the University of Southern California Rossier School of Education), David Drew (then dean of the School of Educational Studies at Claremont Graduate University), and Darryl Yong. Pam Mason became the Executive Director of MfA LA in January 2008.

Since then, MfA LA has supported 179 secondary school teachers in the greater Los Angeles area through its two multi-year teacher fellowship programs. It currently supports 17 Early Career Fellows and 70 Master Teacher Fellows. Like most other MfA sites, MfA LA began by emulating MfA's original Teacher Fellowship program, a five-year fellowship for beginning secondary-school mathematics teachers that included one year of teacher credentialing. MfA LA Teacher Fellows could choose to become credentialed at either the University of Southern California or Claremont Graduate University. In 2013, MfA LA began awarding five-year Master Teacher Fellowships for experienced teachers. In 2014, the Teaching Fellowship program was replaced by a four-year Early Career Teaching Fellowship, which was identical in every way, except that it did not include support for teachers to become credentialed. (For the purposes of this article, we refer to all beginning teachers currently supported by MfA LA as Early Career Teaching Fellows.) Given that mathematics teachers were responding to increasing demands for secondary-school CS teachers, MfA LA began offering Master Teacher Fellowships to CS teachers in 2018.

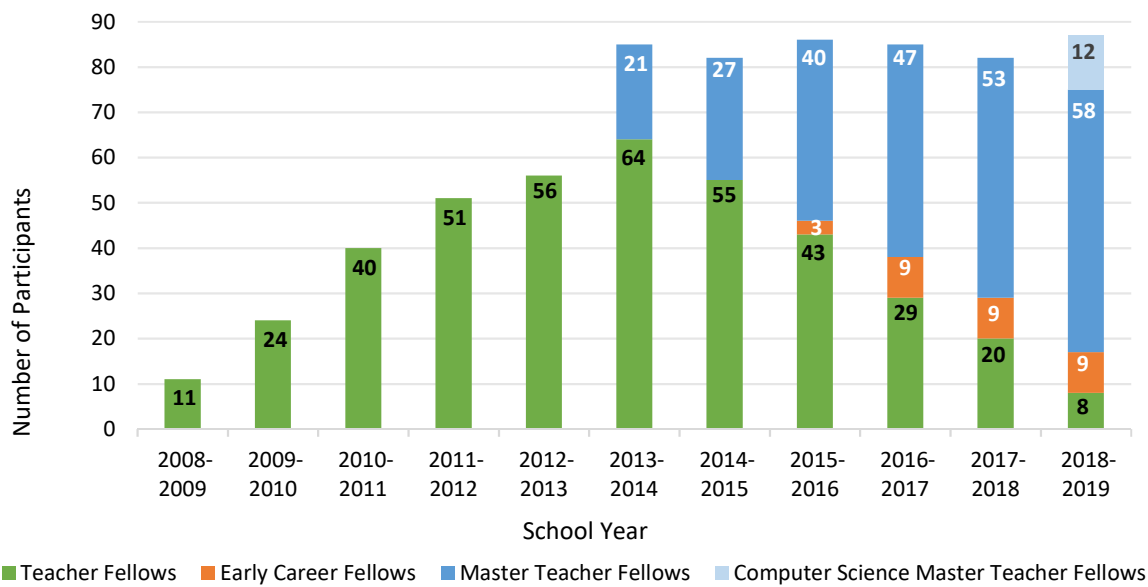

Figure 1. MfA LA History and Number of Participating Teachers

MfA LA receives support from four grants from the National Science Foundation Robert Noyce Teacher Scholarship Program (DUE 0934923, 1136415, 1439917, 1758455), Math for America, philanthropic organizations, and individual donors.

## Characteristics of Successful Professional Development Programs

During the early days of MfA LA, we (Pam and Darryl) primarily drew on our prior experiences while designing its programs. Over the years, we began to recognize how those design choices reflected wisdom from the rich literature about teacher professional development and sought to bring the former in line with the latter.

Over the last 30 years, research on the characteristics of effective professional development, comprised of research commentaries, empirical studies, and meta-analyses, has been remarkably self-consistent (Banilower, Boyd, Pasley, & Weiss, 2006; Birman, Desimone, Porter, & Garet, 2000; Blank & de las Alas, 2009; Borko, Jacobs, & Koellner, 2010; Darling-Hammond, Hyler, & Gardner, 2017; Desimone, Porter, Garet, Yoon, & Birman, 2002; Garet, Birman, Porter, Desimone, & Herman, 1999; Garet, Porter, Desimone, Birman, & Yoon, 2001; Little, 1993; Loucks-Horsley et al., 1987; Loucks-Horsley, Stiles, Mundry, Love, & Hewson, 2010; Stein, Smith, & Silver, 1999; Timperley, 2008; Timperley & Alton-Lee, 2008; Wilson, 2013). Each of these citations contains a set of characteristics of effective professional development. We reproduce the oldest and most recent of these cited below to illustrate the consonance between them, save for some variation in the content foci of professional development efforts.

| | |
|---|---|
| A1. Collegiality and collaboration | B1. Is content focused |
| A2. Experimentation and risk taking | B2. Incorporate active learning |
| A3. Incorporation of available knowledge bases | B3. Supports collaboration |
| A4. Appropriate participant involvement in goal setting, implementation, evaluation, and decision making | B4. Uses models of effective practice |
| | B5. Provides coaching and expert support |
| A5. Time to work on staff development and assimilate new learnings | B6. Offers feedback and reflection |
| | B7. Is of sustained duration |
| A6. Leadership and sustained administrative support | |
| A7. Appropriate incentives and rewards | |
| A8. Designs built on principles of adult learning and the change process | |
| A9. Integration of individual goals with school and district goals | |
| A10. Formal placement of the program within the philosophy and organizational structure of the school and district | |

Table 1. Characteristics of effective professional development, according to Loucks-Horsley et al., 1987 (left) and Darling-Hammond, Hyler, and Gardner, 2017 (right).

As we describe the specific features of MfA LA in the next section, we refer back to these characteristics using the letters and numbers in Table 1 to illustrate how they are reflected the design of MfA LA's programs. It should be noted, however, that these characteristics are necessary but not sufficient conditions for effective professional development (Kennedy, 2016; Sztajn, Borko, & Smith, 2017).

## 3. MfA LA's Program Design

Though schools and districts are the primary provider of professional development for in-service teachers in the United States, there is a rich history of other types of organizations supporting teachers' professional learning. In her 1993 paper, Judith Warren Little presents a taxonomy of professional development organizations not run by schools and districts: (1) subject collaboratives and other networks, (2) subject matter associations like NCSM and NCTM, (3) school-university collaborations targeted at school reforms, and (4) special institutes and centers. Based on this taxonomy, we classify MfA LA as a school-university, subject-specific teacher collaborative.

In this section, we describe the design of MfA LA's programs following recommended practices for reporting on teacher professional development programs (Loucks-Horsley et al., 2010; Sztajn, 2011).

### Context
Like many other dense urban areas, the greater Los Angeles region is a challenging place in which to be a public school teacher. Los Angeles frequently appears on Forbes' lists of America's Most Overpriced Places. The high cost of living in Los Angeles creates financial challenges for teachers who want to live near where they teach. Additionally, Los Angeles is also highly racially and economically segregated (Silver, 2015; Stokes, 2018). According to 2018 U.S. Census Bureau estimates, roughly 60% of people (above the age of 5) in Los Angeles County speak a language other than English at home.

During the 2018–19 academic year, MfA LA teachers taught in 60 different schools and 21 different school districts in the greater Los Angeles area. Almost 60% of them work in Los Angeles Unified School District (LAUSD). Second in size only to the New York City School District, LAUSD consists of 1,322 schools (including 216 charter schools) that serve nearly 700,000 students. Other districts in this area tend to be quite small and serve specific local communities. About 8% of MfA LA teachers work in public charter schools, the rest in traditional, non-charter public schools. Even though MfA LA teachers work in wide range of school contexts, they all face a common challenge of teaching mathematics and CS in urban, public school settings. The average school of MfA LA Fellows has high proportions of Hispanic/Latin[1] students (68%) and students eligible for free or reduced lunch (74%).

Like many other places around the country, California is facing a severe teacher shortage (Darling-Hammond, Sutcher, & Carver-Thomas, 2018). Teacher attrition accounts for 88% of the demand of new teachers in California (Darling-Hammond et al., 2018). This shortage compounds the problem that low-income and historically marginalized students tend to have less qualified teachers (Adamson & Darling-Hammond, 2011). Half of new California math and science teachers are entering the profession without an undergraduate degree in their respective fields (Darling-Hammond et al., 2018; Ingersoll, 2002). In contrast, only 21% of MfA LA teachers do not have an undergraduate degree in their content area and all of them receive monthly training within their content areas. Teachers, especially those who are our best instructors, need rewards, opportunities for growth, and more autonomy to feel vital and remain in the profession (Fiarman, 2007; Ingersoll & May, 2010; Kirkpatrick, 2007). MfA LA addresses these issues by helping to retain talented, trained teachers in the profession and to support experienced teachers in mentoring beginning teachers into the profession.

### Knowledge and Beliefs about Teacher Learning

MfA LA adopts sociocultural perspectives on teachers' learning (Cobb, 1994; Lave & Wenger, 1991; Stein et al., 1999) because of the relevance of these ideas to teachers' communities of practice. Instead of only considering how participation in MfA LA's programs affect teachers' skills, knowledge, and behaviors, we are also interested in the development of shared conceptions of teaching, students, and content knowledge that arise within the organization. We attend to these sociocultural factors because we believe that MfA LA Fellows learn about teaching through participation in communities of practice with other teachers and that learning is shaped by and shapes their social relations with other Fellows. Stories that teachers share over lunch with each other about their classrooms can be as impactful on their teaching practices as the transmission and assimilation of new instructional strategies or content knowledge.

What are the implications of a sociocultural perspective of teachers' learning on MfA LA's programs? If we posit that knowledge resides in the actions of people in a community of practice and evolves as individuals participate in that community, then the conditions necessary for fruitful participation in community become our top concern. We gather together teachers who have a genuine desire to advance their own teaching practices, enculturate them into a community of practitioners who openly share their ideas and resources with each other, allow them the chance to shape their own learning, and reduce the barriers to participation in this community as much as we can.

---

[1] We use the term "Latin" here as a gender-neutral alternative to "Latino".

## Critical Issues

According to Sztajn, critical issues faced by teacher professional development programs should also be reported when discussing the results of these efforts. These critical issues include the program's content focus (Section 1), degree to which teachers from same school participate, information about providers (Section 2), use of incentives (Section 3), compulsory or voluntary participation by teachers (Section 4), and teachers' voice in decision making about their own learning professional development (Sztajn, 2011, p. 229). Some of these have been reported elsewhere (as indicated); the remainder will be discussed in this section.

About 70% of MfA LA teachers teach at a school with at least one other current MfA LA Fellow or program alum. That is both the result of informal word-of-mouth communication (teachers telling their colleagues about open positions or principals calling to inquire about hiring teachers) and also a deliberate design of our Master Teacher Fellowship, which requires teachers to apply together to the program. This "clustering" of teachers at school sites is an important strategy for increasing teachers' job satisfaction and effectiveness, particularly because MfA LA teachers are geographically spread out over an enormous area that makes it difficult for them to otherwise meet socially or professionally.

Because teachers generally have few opportunities to shape the kind of professional development that they receive from schools and districts (Little, 1993), they often complain about professional development that is "done to them" as opposed to "designed for them." In contrast, MfA LA programs and professional development meetings are designed based on teachers' needs and concerns, which are regularly solicited via anonymous feedback surveys, conversations between program staff and Fellows, and current national and regional trends in the teaching of mathematics and CS. Consequently, unlike many other professional develop programs, MfA LA does not restrict its professional development to particular instructional strategies or curricula. However, three common themes in MfA LA's professional learning are (1) content knowledge for teaching (Loewenberg Ball, Thames, & Phelps, 2008), (2) high-quality instruction in mathematics and CS (K-12 Computer Science Framework, 2016; National Council of Teachers of Mathematics, 2014), and (3) equity, diversity, and inclusivity in the classroom (Bartell et al., 2017; Goffney & Gutiérrez, 2018; Jackson & Cobb, 2010).

## Theory of Action

MfA LA's programs rely a theory of action that providing the right kinds of professional supports for teachers will help them become *teacher leaders*, which the National Board for Professional Teaching Standards defines "instructional and organizational change agents who have a critical impact on school, teacher and student success" (2012). The goal is for these teacher leaders to not only be highly-effective teachers, but to also shape school policy, define outcomes and standards, choose and implement curricula, design professional learning for their colleagues, invent and share teaching innovations, and take on numerous other crucial roles in the educational ecosystem (Barth, 2001; York-Barr & Duke, 2004).

| Inputs<br>What is invested in MfA LA | Outputs | | Outcomes | | |
|---|---|---|---|---|---|
| | Activities<br>What MfA LA does | Participation<br>Who MfA LA reaches | Short-term | Medium-term | Long-term |
| Program Staff<br><br>Higher-Ed Faculty Associates<br><br>Research Partners<br><br>Steering Committee<br><br>Advisory Board<br><br>Financial Supports<br><br>Time | Fellows receive stipends + collab. periods (MTFs only) | Fellows participate in professional support and learning | Fellows deepen content knowledge | Student learning attitudes and outcomes improve | Lasting instructional improvements at schools |
| | Fellows participate in PD activities & conferences & NBC | Students learn and participate in fellows' classrooms | Fellows' repertoire of instructional strategies expands | Fellows' pedagogical judgment and reasoning sharpen | Fellows recognized for their expertise and contributions |
| | Fellows teach in classrooms and work with their colleagues | Schools/districts support fellows and benefit from MTF improvement projects | Fellows strengthen advise-seeking & support networks | Fellows positively influence their depts. and schools | Fellows retained in profession |
| | Fellows sharing expertise with other fellows | Teaching community receives leaders and new ideas | Fellows improve sense of wellness | Fellows enhance their leadership capabilities | Fellows share expertise and ideas outside of MfA LA |
| | Coaches support teachers in their classrooms | Research partners collaborate with MfA LA on educational research | Fellows feel respected, supported and validated as teachers | MTFs receive National Board Certification | MfA LA helps to contribute to knowledge base on effective teacher |
| | MTFs implement improvement projects in schools | | | Fellows' confidence and job satisfaction increase | |

Figure 2. MfA LA's logical model clarifies the design of its programs and evaluation (MTF stands for Master Teacher Fellows; NBC stands for National Board Certification; PD stands for professional development)

## Programmatic Details

In this section, we describe features of MfA LA's programs and connect them to characteristics of effective professional development from Table 1, indicated using parentheses.

MfA LA currently offers two types of fellowships: Early Career Teacher Fellowships and Master Teacher Fellowships, for beginning and experienced teachers, respectively. The Early Career Teaching Fellowship is a four-year commitment and the Master Teacher Fellowship program is a five-year commitment (B7). Both programs are designed to give secondary school mathematics and CS teachers the training, resources and connections they need to become leaders who can enact school change, broaden access to high-quality education and raise student performance.

Both fellowship programs offer teachers a $10,000 annual stipend for participating in regularly monthly professional development meetings together and completing other program expectations (A7). These financial incentives follow the recommendation from the National Mathematics Advisory Panel in its 2008 Final Report that teacher pay may be an effective way to attract and retain effective teachers. MfA LA evaluations have consistently revealed that its Fellows do not stay in the program primarily for the financial compensation, but instead for the professional learning and community that they receive. However, Fellows report this annual salary supplement helps them live more balanced lives, for example, by enabling them to live closer to school and commute less, to start a family, or to pay off college loans.

Both programs provide teachers with funding to attend regional and national conferences (A1, B3), for classroom supplies and equipment, and training and support to become certified by the National Board of Professional Teaching Standards, if they have not yet been certified. Conference attendance allows MfA LA fellows to connect with a broader community of educators and provides teachers opportunities to develop their leadership capacity by presenting their work to others (A6).

A key idea of the Master Teacher Fellowship Program is that there is no better way for teachers to develop their leadership capacity than to directly engage in the work of designing and orchestrating student-centered changes in their schools and districts, with the support of other experienced teachers and STEM educators (A1, A2, A4, A9, B3). Applicants must identify some critical need relating to CS and/or math education in their school and propose a plan to address those needs; these plans also need to include some measurable outputs by which success of implementation can be assessed. MfA LA requires that pairs of teachers working at the same school site apply to the program and, if selected, commit to participate in the program together for all five years of the fellowship so that they will have a greater chance of bringing about real, sustained change in their departments.

Time is a crucial, scarce resource for teachers. Expecting them to take on additional creative and instructional responsibilities in addition to all that they are already doing is a recipe for fatigue and burnout (Barth, 2001). To enable Master Teachers to implement their plans for school site improvement, MfA LA coordinates with each school's administrators to buy out a period for partnering teachers so that they can have a collaboration period in addition to their own "conference" periods (A1, A5, B3). A common collaborative period is only deployed if it is in the best interest of the partnering teachers and their school. Coordination with school administrators is also a mechanism to ensure that teachers' goals align with school and district goals (A9).

Monthly all-day Saturday professional development meetings are the primary mechanism by which MfA LA teachers learn from one another and establish a close-knit community (A1, A8, B1, B2, B3). These meetings are interactive, tightly-orchestrated meetings in which teachers work together in varied configurations, based on their interests, needs, and school contexts. These monthly meetings have evolved over the last 10 years, but have largely retained the same structure: an opening activity designed to help teachers deepen their disciplinary content knowledge, workshops designed to help teachers deepen their own teaching practice, and time and structures for teachers to work and share with each other. The opening content area activities are a highlight for many MfA LA Fellows; these activities often do not directly connect to their school subject matter content, but instead are designed to engage teachers in doing mathematics and/or CS together. STEM faculty from partnering institutions participate in the design and delivery of these activities. The specific topics included in these whole-day meetings depend in large part on based on MfA LA Fellows' needs and interests (A4).

During these monthly meetings, all teachers participate in a year-long Working Group, which function somewhat like a professional learning community (Stoll, Bolam, McMahon, Wallace, & Thomas, 2006), except that that members of the group may not all work at the same school site (A1, A2, A5, B3, B7). In these Working Group meetings, teachers often share classrooms artifacts such as student work samples, video recordings, or lesson plans with each other. They plan together, troubleshoot issues, and share resources with each other. Fellows who attend conferences are expected to implement any new ideas acquired and share that with the rest of the MfA LA community (A1, A2, B3).

MfA LA also offers coaching support to provide teachers with feedback on their classroom practices (A8, B5, B6). As this feedback is more proximal to teachers' instructional practices, it complements the Working Groups, workshops during monthly meetings, and conference attendance, which are increasingly distal, though still useful, to teachers' own classroom contexts. Veteran teachers hired by MfA LA provide coaching support. MfA LA also partners with researchers from Vanderbilt University (through NSF HER-1620920) to offer a small subgroup of Master Teachers coaching through an innovative video-based formative feedback (Jurow, Horn, & Philip, 2019).

Finally, MfA LA has actively cultivated an intimate, caring community of teachers in a variety of ways. For example, MfA LA retains the services of a therapist to provide counseling services to its Fellows. In 2017, MfA LA began setting aside a time during monthly meetings for teachers to work in small groups around wellness goals. A warm, affectionate spirit is exemplified in the way that teachers are aware of and celebrate each other's personal and professional accomplishments, the amount of time they spend with each other outside of formal MfA LA events, and the support that they offer each other.

Selection of new MfA LA Fellows takes place each spring. After an initial review of applications, finalists undergo interviews and classroom observations. Selection criteria include passion and aptitude for mathematics and/or CS; demonstrated effectiveness at teaching; desire to advance their own teaching practices; understanding of diversity issues and desire to broaden participation in math and CS; ability to work with and learn from others; leadership capacity and sense of personal responsibility; and ability to adapt to changing situations and unexpected events. Master Teacher applicants must also demonstrate potential for enacting school change through their proposed improvement plan. These plans are judged on their appropriateness, potential for broadening access to quality mathematics and CS education, alignment with school and district standards, thoughtfulness of design, and potential for sustainable change.

MfA LA has a good record of assembling a broadly representative set of Fellows, but continues to actively encourage teachers of color (particularly, African-American and Hispanic/Latinx) to apply. At present, 57% of all MfA LA Teacher Fellows identify as female, 5% as African-American, 26% Hispanic/Latin, and 5% as multiracial.

## 4. Evidence of Impact

Before describing what we know about the impact of MfA LA's programs on teachers and students, it is important to highlight that the incentives and mechanisms for teachers' participation in the professional learning offered by MfA LA are different than that offered by schools and districts. Often, teachers are required to attend school- or district-run professional development. In contrast, MfA LA is a non-profit organization that offers competitive multi-year fellowships to teachers. MfA LA fellows opt in and are selected by the organization. This mutual selection between MfA LA and its Fellows creates a unique context for professional learning.

Cobblestone Applied Research & Evaluation, Inc. conducts on-going and comprehensive evaluation of MfA LA's programs. A summary of all sources of data for research and evaluation appears in Table 2. Using MfA LA's logical model (Figure 2), MfA LA worked with Cobblestone to define 13 evaluation questions and map data sources onto each question.

| **Data involving Fellows' students and non-Fellows' students:** |
|---|

> - Annual survey of student math attitudes at beginning and end of each school year †
> - California Standards Tests (CST) student achievement data (2010 to 2015) †
> - California High School Exit Examination pass rates (2010 to 2015) †
> - Smarter Balanced Assessment System (SBAC) student achievement data (2016 and onwards) †
>
> **Data involving teachers:**
> - Annual surveys of principals of Fellows †‡
> - Annual surveys of Fellows †‡
> - Program pre- and posttest surveys of Fellows †‡
> - Master Teacher work logs (every semester) †‡
> - Classroom observations conducted by Cobblestone (at beginning, middle, end of fellowship) †‡
> - Employment status of current and past Fellows ‡
> - Structured interviews of Fellows as they begin and finish their fellowship ‡
> - Annual structured interviews of Master Teachers about their projects ‡
>
> **Data involving the program as a whole:**
> - Monthly and annual surveys administered to Fellows about the program †‡
> - Observations of all monthly professional development meetings by Cobblestone ‡

Table 2. Sources of data for MfA LA's program evaluation (†= quantitative data, ‡=qualitative data)

## Teacher Effects

MfA LA has strong evidence that the support it provides teachers helps to keep them from leaving the teaching profession. Across LAUSD secondary schools, 61% of new teachers are still teaching in LAUSD three years later; across the United States, 56% of new teachers are still teaching five years later (Ingersoll, Merrill, Stuckey, & Collins, 2018). In contrast, of the new teachers that joined MfA LA from 2008 to 2013, 90% of them were still teaching secondary school three years after they started and 66% of them were still teaching in 2018. These numbers are all the more significant considering that attrition rates for mathematics teachers are higher than for many other disciplines and that attrition in high-need schools is also higher. (See Section 3 for a general description of schools that MfA LA Fellows teach in.)

Fellows report to us that the professional and financial support they receive are the two main reasons why MfA LA helps them to stay in the profession. One teacher reported, "… being a teacher without the support of MfA LA would have been very challenging… I don't think I would've been able to survive more than two years of teaching if I didn't have a community of educators." Another teacher described the effect of MfA LA's financial support: "The MfA LA stipend has become my de facto annual savings; without it I would be living paycheck to paycheck, and my wife would have to take out more in student loans. Without the stipend, I would be compelled to seek alternative income sources, and a relocation or career change would not be ruled out."

There is growing evidence that MfA LA Master Teacher Fellowship Program is helping Master Teachers increase their leadership skills. Nearly two-thirds of MfA LA Master Teachers are certified by the National Board for Professional Teaching Standards; the vast majority of them acquired that certification with support from MfA LA. During the 2017-18 school year, on average each Master Teacher attended 139 hours of professional development, spent 49 hours preparing for professional development, and another 44 hours leading professional development.

Cobblestone includes two survey instruments measuring transformational leadership (Podsakoff, MacKenzie, Moorman, & Fetter, 1990) and instructional leadership (Hallinger & Heck, 2011) for Master Teachers at the beginning and end of their five-year fellowships. Due to small sample size, data from 2012-2014 Master Teacher cohorts was combined and a Wilcoxon Signed-Ranks tests were conducted on both scales. The analysis used list wise deletion, so that only teachers who completed both a pretest and posttest were included. Overall, the differences between the average pretest and posttest ratings for transformational leadership were found to be statistically significant ($p \leq 0.01$), whereas the instructional leadership ranks were not ($p > 0.05$).

| Construct | Pretest Mean *(SD)* | Posttest Mean *(SD)* |
|---|---|---|
| Transformational Leadership | 3.59 *(0.26)* | 4.16 *(0.40)* *** |
| Instructional Leadership | 3.60 *(0.47)* | 3.76 (0.60) |

Table 3. Overall changes in reported transformational and instructional leadership over five years by the 2012, 2013, and 2014 Master Teacher cohorts (N=17). Respondents answered questions on a 5-point Likert scale, where 1 = Strongly Disagree and 5 = Strongly Agree and 1 = Never and 5 = Always, respectively.

Furthermore, MfA LA seems to be helping Teacher Fellows have an impact on their colleagues and schools. Anecdotally, about two-thirds of all Fellows reported taking on some form of leadership at their schools during the 2018-19 school year: serving on a school-governance committee, being a department chair, providing professional development for colleagues at their schools, and so on. Here are some representative quotes from teachers:

> "The opportunity to give presentations [at MfA LA] regularly has helped me to become a better and more confident public speaker."
>
> "I have confidence in my teaching practice and pedagogy and so I feel able to support other teachers at my school."
>
> "Being in a partnership with [another Fellow] has helped me to develop my leadership skills. We push each other to lead PDs at school and be supportive of other teachers. We are able to do this because we have the time and support provided by MfA LA."
>
> "Participating in MfA LA has had a very positive impact as a mathematics leader, especially at my school and even more this year with the district. I have established a reputation where others value and seek my input. After joining MfA LA and doing National Board Certification, I became much more confident and vocal about mathematics education."
>
> "I became department chair for the first time in my school this year. Also I became a part of the School Site Committee and the Instructional Team of the school. I would never have even considered these leadership roles if not for MfA LA."

Cobblestone also independently observes Fellows' classrooms at the beginning, middle, and end of each person's five-year fellowship using a modified version of a classroom observation tool by Tisdel and Ehlert (1998). Our goal is not to call out particular teachers or evaluate their teaching, but to describe whether cohorts of teachers as a whole are making progress in their instructional practices. Table 4 shows an increase in all subscales of the observation protocol. The data are aggregated for the 2010-13 Early Career Fellows and 2012-13 Master Teacher Fellows at the beginning and end of their fellowships.

| Subscale | EC Pre Mean | EC Post Mean | MT Pre mean | MT Post Mean |
|---|---|---|---|---|
| Mathematics Content | 4.29 | 4.55 | 4.28 | 4.72 |
| Pedagogy | 3.61 | 4.41 | 3.99 | 4.71 |
| Classroom Management | 3.75 | 4.48 | 4.39 | 4.62 |
| Classroom Culture | 3.68 | 4.43 | 4.08 | 4.75 |
| Physical Environment | 4.21 | 4.71 | 4.75 | 5.00 |

Table 4. Mean classroom observation subscale scores for Early Career (EC, n=32) and Master Teacher (MT, n=7) Fellows at the beginning and end of their fellowship periods.

Principals also quite uniformly reported to us that MfA LA Fellows made a positive impact on their school. Table 5 summarizes the mean scores reported by principals who completed surveys in 2018.

| Please rate the extent to which you agree with the following statements about the M*f*A LA program. | Mean (*SD*) |
|---|---|
| I would consider hiring an additional M*f*A LA teacher in the future | 4.52 (*0.75*) |
| Teachers' participation in M*f*A LA has had a positive impact on students at my school | 4.52 (*0.81*) |
| M*f*A LA teacher(s) at my school site has had a positive impact on the math department curriculum | 4.50 (*0.96*) |
| M*f*A LA teacher(s) at my school site have positively influenced the pedagogy used in the math department | 4.41 (*0.85*) |
| Teachers' participation in M*f*A LA has had a positive impact on teachers at my school | 4.41 (*0.96*) |
| [If you have any M*f*A LA Master Teacher Fellows at your school site], having a common planning period for Master Teacher Fellows is important for the success of their work together. | 4.40 (*0.52*) |
| It has been challenging having an M*f*A LA teacher at my school. | 2.32 (*1.55*) |

Table 5. Mean responses of principals (N=22) about MfA LA Fellows at their schools; Scale: 1 = Strongly Disagree to 5 = Strongly Agree.

## Student Effects

Given the connection, though weak, between students' attitudes toward mathematics and their achievement in mathematics (Ma & Kishor, 1997), we sought to measure MfA LA Fellows' effectiveness at improving students' attitudes toward learning mathematics. Each year, we select a random sample of MfA LA Fellows and ask their school administrators for consent to collect student attitude data. Students of MfA LA Fellows and students of other teachers at their schools teaching similar courses were then asked to consent to collection of attitudinal data at the beginning and end of the school year. We use the Minnesota Mathematics Attitude Inventory (Sandman, 1980), an instrument we found to have good reliability at both pre- and posttest.

For brevity, we report most recent results from the 2017-18 academic year. During that year, our sample included 744 students of MfA LA Fellows (treatment) and 807 students of non-MfA LA teachers (control). Chi-square tests revealed that treatment and control students were equivalent on sex and ethnicity. Multivariate analysis of covariance (MANCOVA) was used to compare the posttest mathematics attitudes of treatment and control students, controlling for students' attitudes as pretest as well as demographic characteristics. This analysis revealed that students of MfA LA Fellows reported significantly higher perceptions of their teacher, more positive attitudes about their enjoyment of mathematics and value of mathematics at posttest. However, effect sizes were negligible. Overall results of the posttest survey for each of the six key dimensions of this inventory are shown in Table 6. These findings are similar to those from previous years.

| Subscale | MfA LA students (n=744) | Control students (n=807) |
|---|---|---|
| Math Teacher | **3.20** | **2.99** |
| Math Anxiety* | 2.87 | 2.82 |
| Math Value | **2.95** | **2.87** |
| Math Self-Concept | 2.68 | 2.62 |
| Math Enjoyment | **2.58** | **2.46** |
| Math Motivation | 2.08 | 2.01 |

Table 6. Mean student attitudes reported on six subscales of the Minnesota Mathematics Attitude Inventory, at the end of the school year. Bolded values indicate a significant difference (p<0.05); 1=Strongly disagree, 4=strongly agree; * math anxiety was reverse coded so higher scores indicate less anxiety toward mathematics.

While improving student standardized achievement scores is not a primary goal of MfA LA, the way we currently operationalize student learning is through these scores in the absence of other, more valid data sources. Therefore, we have investigated how students of MfA LA Fellows perform relative to other students on standardized tests. Additionally, standardized measures for learning in CS do not yet exist, so we only report standardized mathematics achievement scores. Prior to 2015, MfA LA collected California Standards Tests (CST) scores from districts in which our Fellows taught. After 2016, we had access to Smarter Balanced Assessment System (SBAC) student achievement data from districts.

For brevity, we report most recent SBAC results from the 2017-18 academic year. Anonymized student data were received from 10 middle schools and 29 high schools representing 4 districts. Out of an initial sample of 13,817 students, propensity score matching (Rosenbaum & Rubin, 1985; Thoemmes, 2012) was used to generate two subsamples of students of MfA LA Fellows (treatment) and non-MfA LA Fellows (control), matched on demographic markers (such as gender, ethnicity, and English proficiency) and school-level variables (such as percentage of English language learners and students receiving free or reduced lunch). This matching produced 781 middle school and 1,092 high school students in both treatment and control groups. It should be noted that though there were more high school students in the sample, SBAC data is only available for Grades 6-8 and 11. Consequently, the sample of high school students is a somewhat skewed proportion of students since only a fraction of MfA LA high school Fellows teach 11th graders.

|  | MfA LA Students | | Control Students | | F (df) | $\eta^2$ |
|---|---|---|---|---|---|---|
|  | n | M (SD) | n | M (SD) |  |  |
| **Middle School** | 781 | 2590.13 (107.78) | 781 | 2500.72 (112.71) | 187.53 (1, 1556)*** | .11 |
| **High School** | 1,092 | 2593.28 (114.36) | 1,092 | 2579.55 (119.91) | 11.44 (1, 2178)** | .01 |

Table 7. One-way ANOVA results for mean comparisons of SBAC scores; **p<.01; ***p<.001; Partial $\eta^2$ is a measure of effect size that represents the proportion of variance in the outcome that is accounted for by the predictor (whether a student has a MfA LA teacher or not).

|  |  | Middle School (*n*=1,562) | High School (*n*=2,184) |
|---|---|---|---|

| Level | Description | MfA % (n) | Control % (n) | MfA % (n) | Control % (n) |
|---|---|---|---|---|---|
| 4 | Exceeded Standard | **38.9% (303)** | **15.2% (119)** | 14.8% (162) | 11.7% (128) |
| 3 | Standard Met | 20.2% (158) | 15.7% (122) | 23.3% (254) | 25.1% (274) |
| 2 | Standard Nearly Met | 24.2% (189) | 26.5% (207) | 27.6% (301) | 26.6% (291) |
| 1 | Standard Not Met | **16.7% (131)** | **42.6% (333)** | 34.3% (375) | 36.5% (399) |

Table 8. SBAC achievement level rates of MfA LA students compared to control students. Statistically significant differences in achievement level rates with p<.05 are bolded.

Table 7 presents mean SBAC scores whereas Table 8 presents the proportion of students in different achievement levels on the SBAC. As can be seen in Table 7, both middle and high school students of MfA LA Fellows scored higher on the SBAC than control students and these differences were statistically significant: Middle school $F(1, 1556) = 187.5$, $p < .001$; High school $F(1, 2178) = 11.4$, $p < .01$.

The figures in Table 8 indicate that a higher proportion of MfA LA students than control students achieved Level 4 (38.9% vs 15.2%, respectively), and a lower proportion of MfA LA students than control students at Level 1 (16.7% vs 42.6%, respectively). The effect size for middle school was moderate (Cramer's $V = .33$), suggesting that students who have an MfA LA Fellow as an instructor are somewhat more likely to achieve a higher SBAC level. For high school students, there were no statistically significant differences in the proportion of MfA LA and control students across achievement levels, though there have been significant findings in previous years.

Taken as a whole, these findings help to validate MfA LA's theory of action (Figure 2) and show that these efforts are helping teacher stay in the profession, become more effective teachers, become more effective leaders, and improve students learning outcomes.

## 5. Lessons Learned

This retrospective article might give one the false impression that MfA LA's programs unfolded linearly from core principles and effective professional development practices to logic models and programmatic choices. A more accurate description of MfA LA's progress is more chaotic and includes missteps and growth in fits and starts. However, MfA LA has always grown organically in response to what our Fellows tell us they need, shifting conditions in schools and districts, and constant program evaluation.

Perhaps the most important thing that we have learned is that good professional development requires significant time and resources and that those two things are constantly in tension with each other. Teacher learning takes time, and being able to measure that learning takes even more time. Therefore, any significant professional development effort needs to measure progress over years, not weeks and months. Another example of this tension derives from the reality that teaching in public schools is demanding and time consuming. Buying out collaborative periods is one way to give teachers the time they need to make significant progress on their instructional goals, but it is expensive.

Another tension common in many professional development efforts is the balance between scale and community. MfA LA teachers report that the community and holistic support they receive provide solace from constantly shifting policies and leadership at their schools and districts. And yet, professional development programs are often scrutinized and encouraged to scale up for greater impact. How can effective professional development be scaled up while preserving its integrity and robust sense of community?

And of course, a constant challenge for most professional development efforts is finding and sustaining financial support. A few years ago, MfA LA was just one of six satellite MfA programs created outside of New York. However, besides Los Angeles and New York, at present all other MfA sites have discontinued or are winding down their programs. MfA LA has remained viable due to its fundraising successes. That fundraising was made possible in part by the robust program evaluation, which allowed MfA LA to communicate its progress and successes with funders and stakeholders. MfA LA now faces the same challenges that many other start-up organizations face as they mature: the problem of convincing funders to contribute to more stable growth instead of more rapid growth, the difficulty of sustaining a coherent mission and strategic vision when founding members of the organization move on, and so on.

In her 1993 paper, Judith Warren Little proposed that "investments [in professional development] beyond the ordinary… are more defensible if they can meet one of these three criteria: (1) they can be credibly tied to a ripple effect (so that the per teacher cost is demonstrably lower than per participant cost); (2) one can claim that the direct individual benefit of this specific program is far more certain than the benefit linked to conventional funding; or (3) the program contributes in demonstrable ways to increased organizational capacity in ways that transcend the impact on those individuals who participate directly in the 'program.'" Time will tell if MfA LA can justify itself in one of these ways.

# References

Adamson, F., & Darling-Hammond, L. (2011). Addressing the inequitable distribution of teachers: What it will take to get qualified, effective teachers in all communities. *Research Brief: Stanford Center for Opportunity Policy in Education*. Retrieved from https://edpolicy.stanford.edu/sites/default/files/events/materials/scope-teacher-salary-brief.pdf

Banilower, E. R., Boyd, S. E., Pasley, J. D., & Weiss, I. R. (2006). *Lessons from a Decade of Mathematics and Science Reform: A Capstone Report for the Local Systemic Change through Teacher Enhancement Initiative.* Retrieved from Horizon Research, Inc. website: http://www.pdmathsci.net/reports/capstone.pdf

Bartell, T., Wager, A., Edwards, A., Battey, D., Foote, M., & Spencer, J. (2017). Toward a framework for research linking equitable teaching with the standards for mathematical practice. *Journal for Research in Mathematics Education*, *48*(1), 7–21.

Barth, R. S. (2001). Teacher Leader. *The Phi Delta Kappan*, *82*(6), 443–449.

Birman, B. F., Desimone, L., Porter, A. C., & Garet, M. S. (2000). Designing Professional Development That Works. *Educational Leadership*, *57*(8), 28–33.

Blank, R. K., & de las Alas, N. (2009). *Effects of teacher professional development on gains in student achievement: How meta-analysis provides evidence useful to education leaders*. Washington, DC: Council of Chief State School Officers.

Borko, H., Jacobs, J., & Koellner, K. (2010). Contemporary approaches to teacher professional development. *International Encyclopedia of Education*, *7*(2), 548–556.

Cobb, P. (1994). Where is the mind? Constructivist and sociocultural perspectives on mathematical development. *Educational Researcher*, *23*(7), 13–20.

Darling-Hammond, L., Hyler, M. E., & Gardner, M. (2017). *Effective Teacher Professional Development* (p. 76). Retrieved from Learning Policy Institute website: https://learningpolicyinstitute.org/product/teacher-prof-dev

Darling-Hammond, L., Sutcher, L., & Carver-Thomas, D. (2018). *Teacher Shortages in California: Status, Sources, and Potential Solutions (research brief)* (p. 61). Retrieved from Learning Policy Institute website: https://gettingdowntofacts.com/sites/default/files/2018-09/GDTFII_Report_Darling-Hammond.pdf

Desimone, L. M., Porter, A. C., Garet, M. S., Yoon, K. S., & Birman, B. F. (2002). Effects of Professional Development on Teachers' Instruction: Results from a Three-Year Longitudinal Study. *Educational Evaluation and Policy Analysis*, *24*(2), 81–112.

Fiarman, S. E. (2007). It's Hard to Go Back: Career Decisions of Second-stage Teacher Leaders. *American Educational Research Association Annual Conference*, 41. Retrieved from http://projectngt.gse.harvard.edu/files/gse-projectngt/files/sef_aera_2007_it_s_hard_to_go_back.pdf

Garet, M. S., Birman, B. F., Porter, A. C., Desimone, L., & Herman, R. (1999). *Designing Effective Professional Development: Lessons from the Eisenhower Program [and] Technical Appendices* (No. ED/OUS99-3). American Institutes for Research.

Garet, M. S., Porter, A. C., Desimone, L., Birman, B. F., & Yoon, K. S. (2001). What Makes Professional Development Effective? Results From a National Sample of Teachers. *American Educational Research Journal*, *38*(4), 915–945. https://doi.org/10.3102/00028312038004915

Garet, M. S., Wayne, A. J., Stancavage, F., Taylor, J., Eaton, M., Walters, K., … Warner, E. (2011). *Middle School Mathematics Professional Development Impact Study: Findings after the Second Year of Implementation. NCEE 2011-4024.* (No. NCEE 2011-4024). Retrieved from National Center for Education Evaluation and Regional Assistance website: https://files.eric.ed.gov/fulltext/ED519922.pdf

Goffney, I., & Gutiérrez, R. (2018). *Rehumanizing mathematics for Black, Indigenous, and Latinx students*. National Council of Teachers of Mathematics.

Guskey, T. R. (2000). *Evaluating Professional Development*. Corwin Press.

Hallinger, P., & Heck, R. H. (2011). Collaborative leadership and school improvement: Understanding the impact on school capacity and student learning. In T. Townsend & J. MacBeath (Eds.), *International Handbook of Leadership for Learning* (pp. 469–485). Springer Netherlands.

Hattie, J. (2009). *Visible Learning: A Synthesis of Over 800 Meta-Analyses Relating to Achievement*. Routledge.

Haycock, K. (1998). *Good Teaching Matters: How Well-Qualified Teachers Can Close the Gap*. Retrieved from http://www.eric.ed.gov/ERICWebPortal/contentdelivery/servlet/ERICServlet?accno=ED457260

Hill, H. (2015). *Review of the mirage: Confronting the hard truth about our quest for teacher development*. Retrieved from National Education Policy Center website: https://nepc.colorado.edu/thinktank/review-tntp-mirage

Ingersoll, R. M. (2002). *Out-of-field teaching, educational inequality, and the organization of schools: An exploratory analysis* (No. R-02-1). Center for the Study of Teaching and Policy.

Ingersoll, R. M., & May, H. (2010). *The Magnitude, Destinations, and Determinants of Mathematics and Science Teacher Turnover*. The Consortium for Policy Research in Education.

Ingersoll, R. M., Merrill, E., Stuckey, D., & Collins, G. (2018). *Seven Trends: The Transformation of the Teaching Force–Updated October 2018*. Retrieved from Consortium for Policy Research in Education website: https://repository.upenn.edu/cpre_researchreports/108

Jackson, K., & Cobb, P. (2010). Refining a vision of ambitious mathematics instruction to address issues of equity. *American Educational Research Association Annual Conference*. Presented at the American Educational Research Association, Denver, CO. Retrieved from http://cadrek12.org/sites/default/files/Jackson%20%20Cobb%20Ambitious%20and%20equitable_vision_100702.pdf

Jacob, A., & McGovern, K. (2015). *The Mirage: Confronting the Hard Truth about Our Quest for Teacher Development.* Retrieved from TNTP website: https://tntp.org/publications/view/the-mirage-confronting-the-truth-about-our-quest-for-teacher-development

Jurow, S., Horn, I. S., & Philip, T. M. (2019). Re-mediating knowledge infrastructures: a site for innovation in teacher education. *Journal of Education for Teaching*, *45*(1), 82–96. https://doi.org/10.1080/02607476.2019.1550607

K-12 Computer Science Framework. (2016). *K-12 Computer Science Framework*. Retrieved from https://k12cs.org/wp-content/uploads/2016/09/K%E2%80%9312-Computer-Science-Framework.pdf

Kennedy, M. M. (2016). How Does Professional Development Improve Teaching? *Review of Educational Research*, *86*(4), 945–980. https://doi.org/10.3102/0034654315626800

Kirkpatrick, C. L. (2007). To invest, coast or idle: Second-stage teachers enact their job engagement. *American Educational Research Association Annual Conference*. Presented at the American Educational Research Association, Chicago, IL. Retrieved from http://projectngt.gse.harvard.edu/files/gse-projectngt/files/clk_aera_2007_paper_to_post.pdf

Lave, J., & Wenger, E. (1991). *Situated learning: Legitimate peripheral participation*. Cambridge university press.

Little, J. W. (1993). Teachers' professional development in a climate of educational reform. *Educational Evaluation and Policy Analysis*, *15*(2), 129–151.

Loewenberg Ball, D., Thames, M. H., & Phelps, G. (2008). Content Knowledge for Teaching: What Makes It Special? *Journal of Teacher Education*, *59*(5), 389–407. https://doi.org/10.1177/0022487108324554

Loucks-Horsley, S., Harding, C. K., Arbuckle, M. A., Murray, L. B., Dubea, C., & Williams, M. A. (1987). *Continuing to Learn: A Guidebook for Teacher Development.* The Regional Laboratory for Educational Improvement of the Northeast and Islands; National Staff Development Council.

Loucks-Horsley, S., Stiles, K. E., Mundry, S., Love, N., & Hewson, P. W. (2010). *Designing Professional Development for Teachers of Science and Mathematics* (Third). https://doi.org/10.4135/9781452219103

Ma, X., & Kishor, N. (1997). Assessing the Relationship between Attitude toward Mathematics and Achievement in Mathematics: A Meta-Analysis. *Journal for Research in Mathematics Education*, *28*(1), 26. https://doi.org/10.2307/749662

Margolis, J. (2010). *Stuck in the shallow end: Education, race, and computing*. MIT Press.

McCaffrey, D., J.R. Lockwood, Koretz, D., & Hamilton, L. (2003). *Evaluating Value-Added Models for Teacher Accountability*. RAND Corporation.

National Council of Teachers of Mathematics. (2014). *Principles to actions: Ensuring mathematical success for all*. NCTM.

Nye, B., Konstantopoulos, S., & Hedges, L. V. (2004). How Large Are Teacher Effects? *Educational Evaluation and Policy Analysis*, *26*(3), 237–257.

OECD. (2005). *Teachers Matter: Attracting, Developing and Retaining Effective Teachers*. OECD.

Owens, A., Reardon, S. F., & Jencks, C. (2016). Income Segregation Between Schools and School Districts. *American Educational Research Journal*, 1159–1197.

Podsakoff, P. M., MacKenzie, S. B., Moorman, R. H., & Fetter, R. (1990). Transformational leader behaviors and their effects on followers' trust in leader, satisfaction, and organizational citizenship behaviors. *The Leadership Quarterly*, *1*(2), 107–142.

Rosenbaum, P. R., & Rubin, D. B. (1985). Constructing a control group using multivariate matched sampling methods that incorporate the propensity score. *The American Statistician*, *39*(1), 33–38.

Sandman, R. S. (1980). The Mathematics Attitude Inventory: Instrument and User's Manual. *Journal for Research in Mathematics Education*, *11*(2), 148–49.

Silver, N. (2015, May 1). The Most Diverse Cities Are Often The Most Segregated. Retrieved June 15, 2019, from FiveThirtyEight website: https://fivethirtyeight.com/features/the-most-diverse-cities-are-often-the-most-segregated/

Stein, M. K., Smith, M. S., & Silver, E. (1999). The development of professional developers: Learning to assist teachers in new settings in new ways. *Harvard Educational Review*, *69*(3), 237–270.

Stokes, K. (2018, July 3). LA's Schools Are Segregated. LAUSD Says There's Only So Much They Can Do. Retrieved June 15, 2019, from LAist website: https://laist.com/2018/07/03/las_schools_are_segregated_lausd_says_theres_only_so_much_they_can_do_about_it.php

Stoll, L., Bolam, R., McMahon, A., Wallace, M., & Thomas, S. (2006). Professional Learning Communities: A Review of the Literature. *Journal of Educational Change*, *7*(4), 221–258. https://doi.org/10.1007/s10833-006-0001-8

Sztajn, P. (2011). Standards for reporting mathematics professional development in research studies. *Journal for Research in Mathematics Education*, *42*(3), 220–236.

Sztajn, P., Borko, H., & Smith, T. M. (2017). Research on mathematics professional development. In J. Cai (Ed.), *Compendium for research in mathematics education* (pp. 213–243). National Council of Teachers of Mathematics.

Thoemmes, F. (2012). Propensity score matching in SPSS. *ArXiv Preprint ArXiv:1201.6385*.

Timperley, H. (2008). *Teacher professional learning and development* (No. Educational Practices-18; p. 32). International Academy of Education.

Timperley, H., & Alton-Lee, A. (2008). Reframing Teacher Professional Learning: An Alternative Policy Approach to Strengthening Valued Outcomes for Diverse Learners. *Review of Research in Education*, *32*(1), 328–369. https://doi.org/10.3102/0091732X07308968

Tisdel, P., & Ehlert, A. (1998). *Mathematics Classroom Observation. Unpublished Instrument.* Retrieved from https://satecsite.org/SATEC/Reports.htm

Wilson, S. M. (2013). Professional Development for Science Teachers. *Science*, *340*(6130), 310–313. https://doi.org/10.1126/science.1230725

York-Barr, J., & Duke, K. (2004). What do we know about teacher leadership? Findings from two decades of scholarship. *Review of Educational Research*, *74*(3), 255–316.